\documentclass[11pt]{amsart}

\newtheorem{lem}{Lemma}[section]
\newtheorem{theo}[lem]{Theorem}
\newtheorem{prop}[lem]{Proposition}

\newtheorem{example}[lem]{Example}
\numberwithin{equation}{section}
\newcommand{\pf}{{\sc Proof:}\quad}
\newcommand{\N}{{\mathbb N}}
\newcommand{\C}{{\mathbb C}}
\newcommand{\loc}{{\mathcal O}}
\newcommand{\fim}{{{\hfill$\rule{2.5mm}{2.5mm}$}}\vspace{1cm}}

\newcommand{\M}{{\mathcal M}}

\newcommand{\ord}{{\rm ord}}

\newcommand{\nd}{\not\hspace{0.1cm}\mid}

\def\N{{\rm I\kern-.18em N}}
\def\GCD{{\rm GCD}}
\usepackage{amsfonts,amssymb,latexsym,indentfirst,eucal,amsmath}
\begin{document}

\pagestyle{myheadings} \markboth{Hefez and Hernandes}{The Analytic
Classification of Plane Branches}

\begin{center}
{\LARGE The Analytic Classification of Plane Branches}
\end{center}
\vspace{5mm}

\begin{center}
{\large  Abramo Hefez}\footnote{partially supported by CNPq,
PRONEX and PROCAD/CAPES.} \ \ \ and  \ \ \ {\large Marcelo E.
Hernandes}\footnote{partially supported by Funda\c c\~ao
Arauc\'aria and PROCAD/CAPES.}

\end{center}
\vspace{1.5cm}

\begin{center}
Abstract
\end{center}

In this paper we give a solution to the open classical problem of
analytic classification of plane branches. \vspace{3mm}

\noindent {\footnotesize AMS Subject Classification (2000).
Primary: 14H20. Secondary: 14Q05, 14Q20, 32S10

} \vspace{1cm}

\section{Introduction}
The aim of this work is to present a solution to the problem of
effective analytic classification of plane branches.

Oscar Zariski, in a course taught at the \'Ecole Polytechnique
\cite{[Z]}, in 1973, inspired by the work of S. Ebey \cite{[E]},
exposed his research on the problem of analytic classification of
plane branches belonging to a given equisingularity class. A great
amount of work is dedicated there to the analysis of particular
examples, showing the need of more effective methods to solve the
problem.

To this respect, in the introduction of \cite{[Z]}, Zariski wrote:

{\em Le probl\`eme de la description compl\`ete de l'\'espace des
modules $M$ d'une classe d'\'equisingularit\'e donne\`e est
enti\`erement ouvert et les quelques exemples du chapitre V
montrent que M a une structure trop complexe pour esp\'erer
r\'epondre totalement \`a la question.

La question, pourtant plus restrictive, de la d\'etermination de
la dimension de la "composante g\'en\'erique" de M n'est pas
r\'esolue.}

The first non-trivial result in this direction was given by C.
Delorme in \cite{[D]}, where he answered the above second question
in a very particular case, describing the generic component of the
moduli space for plane branches with one Puiseux pair and
computing its dimension.

In this paper, we show how one can break the complexity of the
moduli space by stratifying the given equisingularity class by
means of a {\em good numerical invariant} that separates branches
into finitely many types, such that analytic equivalence is each
stratum is manageable.

This is accomplished by considering the sets $\Lambda$ of values
of K\"ahler differentials on branches as finer numerical
invariants than the semigroup of values $\Gamma$ which
characterizes the equisingularity class. If one stratifies with
such sets $\Lambda$ the open set of the affine space representing
the parameter space of an equisingularity class, corresponding to
Puiseux parametrizations with fixed Puiseux pairs, and pick a set
of special representatives in a specified normal form for each
stratum, then the group action that represents analytic
equivalence becomes rather trivial, allowing us to solve in
general the above problems posed by Zariski.

Our setup is similar to that of \cite{[Z]}, adding to it two
techniques with computational flavor. The first one is the use of
a \emph{SAGBI} algorithm, due to L. Robbiano and M. Sweedler
\cite{[RS]}, that we adapted to our situation in \cite{[HH1]} and
\cite{[HH3]}, to describe privileged bases of the local rings of
plane branches as well of the modules of K\"ahler differentials of
these local rings, allowing us to compute the set of invariants
$\Lambda$. The second technique is the algorithm of \emph{Complete
Transversal} due to J. W. Bruce, N. P. Kirk and A. A. du Plessis
\cite{[BKP]} that allows to determine all normal forms of
map-germs under one of Mather's group action, but doesn't allow to
predict a priori what will be the outgoing result. The strength of
our method stems in the conjugation of these two tools that
allows, via the existence of some differentials, to control each
step of the Complete Transversal algorithm, giving explicitly all
possible normal forms and conditions for the analytic equivalence
of germs in normal forms.

All the results we obtain are effective, in the sense that there
is an efficient algorithm that puts any plane branch into its
normal form and it is easy to recognize whether two plane branches
under normal form are equivalent or not. The whole process has
been implemented\footnote{cf. www.dma.uem.br/$\sim$
hernandes/publications.html}. \vspace{5mm}

\section{Preliminaries}

Our ground field is the field $\C$ of complex numbers. Since all
power series we will work with are finitely determined, with
respect to the equivalence relations we will consider, the results
in this work are valid in the formal context and in the analytic
context, as well. In order to have a geometric interpretation of
our objects, we will adopt the analytic point of view. The reader
who desires to find all the known results quoted in this section
is invited to consult \cite{[H]}, where they are gathered with
their proofs.

We denote by $\loc_2$ the ring $\mathbb C\{X,Y\}$ of convergent
power series in two variables with coefficients in $\C$. Let $f$
be in the maximal ideal $\mathcal M_2$ of $\loc_2$ and
irreducible. Then the class of $f$ in $\loc_2$, modulo associates,
is called a {\em plane branch} and denoted by $(f)$. We will
identify $(f)$ with the germ of analytic plane curve at the origin
$\{(x,y)\in (\C^2,0); f(x,y)=0\}$.

We say that two plane branches $(f_1)$ and $(f_2)$ are {\em
equisingular}, writing $(f_1)\equiv (f_2)$ if and only if $(f_1)$
and $(f_2)$ are topologically equivalent as complex immersed germs
of curves in $(\mathbb C^2,0)$; that is, when there exists a
homeomorphism $\Phi:U \to U'$, where $U$ and $U'$ are
neighborhoods of the origin in $\mathbb C^2$ such that $f_1$
(resp. $f_2$) is convergent in $U$ (resp. $U'$) and
$\Phi((f_1)\cap U)=(f_2)\cap U'$. The set of all plane branches
which are equisingular to each other will be called an {\em
equisingularity class}.

When the above transformation $\Phi$ is an analytic isomorphism,
we say that $(f_1)$ and $(f_2)$ are {\em analytically equivalent},
or shortly {\em equivalent}, writing in this case $(f_1)\sim
(f_2)$. So, two branches $(f_1)$ and $(f_2)$ are equivalent if,
and only if, there are an automorphism $\Phi^*$ and a unit $u$,
both of $\loc_2$, such that $ \Phi^* (f_1)=uf_2$.

If we denote by $\loc_f$ the quotient ring $\loc_2/\langle f
\rangle $, where $\langle f \rangle$ represents the ideal
generated by $f$, then one has $(f_1)\sim (f_2)$ if, and only if,
$\loc_{f_1} \simeq \loc_{f_2}$, as $\C$-algebras.

Our main concern in this work is to perform the analytic
classification of plane branches within a given equisingularity
class.

Since the integral closure of $\loc_f$ in its field of fractions
is a complete discrete valuation ring, hence isomorphic to
$\loc_1=\mathbb C\{t\}$, any plane branch $f$ has a
parametrization $\varphi(t)=(x(t),y(t))$ with $x(t)$ and $y(t)$ in
$\M_1$, the maximal ideal of $\loc_1$, not both identically zero,
such that $f(x(t),y(t))=0$. Conversely, any non-zero mapping
$\varphi:(\C,0)\rightarrow (\C^2,0)$, also called a {\em
parametrization}, determines a plane branch. We will call a
parametrization {\em primitive} if it cannot be reparametrized by
a power of a new variable.

It is a well known fact, already used in \cite{[Z]} (see also
Lemma 2.2 in \cite{[BG]}) that, given two plane branches $(f_1)$
and $(f_2)$, parame\-tri\-zed, respectively, by $\varphi_1$ and
$\varphi_2$, then $(f_1)\sim (f_2)$ if, and only if, $\varphi_1$
and $\varphi_2$ are ${\mathcal A}$-equivalent, writing $\varphi_1
\sim_{\mathcal A} \varphi_2$, where $\mathcal A$-equivalence means
that there exist germs of analytic isomorphisms $\sigma$ and
$\rho$ of $(\mathbb C^2,0)$ and $(\mathbb C,0)$, respectively,
such that $\varphi_2=\sigma\circ \varphi_1 \circ \rho^{-1}$.

So, the analytic classification of plane branches reduces to the
$\mathcal A$-classification of parametrizations, which we are
going to undertake in this paper.

Any plane branch is known (cf. \cite{[Z]}) to be equivalent to a
branch with a (primitive) Puiseux parametrization:
\begin{equation}
\varphi(t)=(x(t),y(t))= (t^{\beta_0}, \sum_{i\geq \beta_1}a_it^i),
\end{equation}
with $a_i\in \mathbb C$,  \ $a_{\beta_1}=1$, \ $\beta_0<\beta_1$,
\ $\beta_0\not | \beta_1$ \ and  \ $f\circ \varphi(t)=0$.

The {\em characteristic exponents} of $\varphi$ are the integers
$\beta_0,\beta_1, \ldots$, where for $i\geq 1$,
$$\beta_i=\min\{j; \ a_j\neq 0 \ \hbox{and
$\GCD(\beta_0,\ldots,\beta_{i-1},j)\neq
\GCD(\beta_0,\ldots,\beta_{i-1})$}\}.$$

The associated integers $n_i$ are $n_0=1$ and
$$n_i=\frac{e_{i-1}}{e_i},$$ where
$e_0=\beta_0$ and $e_i=\GCD(\beta_0,\ldots,\beta_i)$.

Since the parametrization is primitive, there is a positive
integer $g$ such that $e_g=1$. This $g$ is called the {\em genus}
of the branch. The Puiseux pairs of $\varphi$ are $(n_i,m_i)$,
$i=1,\ldots, g$, where $\displaystyle m_i=\frac{\beta_i}{e_i}$.

It is a classical result, going back to the thirties, due
essentially to Brauner and Zariski, that the Puiseux pairs form a
complete set of numerical invariants for the topological
classification of plane branches; that is, if $f_1$ and $f_2$ have
parametrizations as (2.1), then $(f_1)\equiv (f_2)$ if, and only
if, $(f_1)$ and $(f_2)$ have the same Puiseux pairs or,
equivalently, the same characteristic exponents.

So, the Puiseux parametrizations as (2.1) with $a_{\beta_i} \neq
0$, for $i=2,\ldots,g$, and such that $a_j=0$, if $\beta_i\leq j
<\beta_{i+1}$ and $e_i\nd j$, determine equisingular branches; and
any plane branch in this equisingularity class is equivalent to
one with such a Puiseux parametrization.

Any parametrization $\varphi$ induces a homomorphism
\[
\begin{array}{cccc}
\varphi^*:& \loc_2     & \to     &\loc_1.\\
        &h    & \mapsto & h\circ \varphi(t)
\end{array}
\]

Assuming $\varphi$ primitive, we define the value $v_\varphi(h)$,
for $h\in \loc_2$, as being $\ord_t(\varphi^*(h))$, the order in
$t$ of the power series $\varphi^*(h)$, defining also
$v_\varphi(0)=\infty$. The semigroup of the extended naturals
$\Gamma_\varphi=v_\varphi(\loc_2)$ will be called the {\em
semigroup of values} of $\varphi$. This semigroup will be
represented by $\Gamma_\varphi=\langle v_0,v_1,\ldots,v_g\rangle$,
where $v_0<\cdots<v_g$ is its minimal set of generators. When
$\varphi$ is a Puiseux parametrization, it is well known (cf.
\cite{[Z]}) that there are relations among these $v_i$'s and the
$\beta_i$'s, given by $\beta_0=v_0$, $\beta_1=v_1$,
$\GCD(v_0,\ldots,v_i)=e_i$, and
\[
v_{i+1}=n_iv_i+\beta_{i+1}-\beta_i.
\]

Hence, the characteristic exponents of  $\varphi$ and
$\Gamma_\varphi$ determine each other.

The semigroup $\Gamma_\varphi$ has a conductor; that is, there is
a natural number $c$ such that $c-1\not \in \Gamma_\varphi$, and
$l\in \Gamma_\varphi$, for all $l\geq c$. This means that the set
$\mathbb N\setminus \Gamma_\varphi$, whose elements are called the
{\em gaps} of $\Gamma_\varphi$, is finite.

So, in order to preserve the Puiseux form given in (2.1), under
the $\mathcal A$-action, the analytic isomorphisms $\sigma$ and
$\rho$ must be of a very special type, as described below.

If $\varphi(t)=(x(t),y(t))$ and
$\varphi_1(t_1)=(x_1(t_1),y_1(t_1))$ are two
para\-me\-tri\-zations as in (2.1), then in order to have
$\varphi_1(t_1)= \sigma \circ \varphi \circ \rho^{-1}(t_1)$, it is
necessary and sufficient that
\begin{equation}
\begin{array}{l}
\sigma(X,Y)=(r^{v_0}X+p,r^{v_1}Y+q),\\
\\
t_1=\rho(t)=rt \sqrt[v_0]{1+\frac{ \varphi^*(p)}{r^{v_0}t^{v_0}}},
\end{array}
\end{equation}
where $r\in \mathbb C^*$ and $p,q\in \loc_2$, with
$v_\varphi(p)>v_0$ and $v_\varphi(q)>v_1$.

The $\mathcal A$-action induced on parametrizations of the form
(2.1) is then given by $(t^{v_0},y(t))\sim_{\mathcal
A}(t_1^{v_0},y_1(t_1))$ if, and only if,
\begin{equation}
y_1(t_1)=r^{v_1}y(\rho^{-1}(t_1))+q(\rho^{-1}(t_1)^{v_0},y(\rho^{-1}(t_1))).
\end{equation}

Ebey and Zariski (cf. \cite{[E]} and \cite{[Z]}) gave some
elimination criteria {\bf (EC)} of parameters of $y(t)$ in a given
parametrization $(t^{v_0},y(t))$ as in (2.1), by means of the
${\mathcal A}$-equivalence.

\quad Let $\varphi=(t^{v_0},t^{v_1}+\sum_{i>v_1}a_it^i)$ be a
Puiseux parametrization, and let $j>v_1$ be an integer. If one of
the following conditions holds, \vspace{1mm}

\noindent {\bf EC1)}\ $j\in \Gamma_\varphi$, \ or \vspace{1mm}

\noindent {\bf EC2)}\ $j+v_0-v_1\in\Gamma_\varphi $, \vspace{1mm}

\noindent then $\varphi$ is ${\mathcal A}$-equivalent to a
parametrization $ (t^{v_0},t^{v_1}+ \sum_{i>v_1}{a'_it^i}),$ with
$a'_i=a_i$, when $i<j$, and $a'_j=0$.

It then follows that any parametrization
$\varphi=(t^{v_0},\sum_{v_1\leq i}a_it^i)$  is ${\mathcal
A}$-equi\-va\-lent to the parametrization $$(t^{v_0},\sum_{v_1\leq
i<c}a_it^i),$$ where $c$ is the conductor of $\Gamma_\varphi$.

Let $\Sigma_\Gamma$ denote the set of all parametrizations of the
above form, such that $\Gamma_\varphi$ is equal to a given
$\Gamma$. This set can be identified with an open set (the
complement of the union of the hyperplanes $a_{\beta_i}=0$,
$i=2,\ldots,g$) of an affine space, whose points are the ordered
sets of the coefficients of $y(t)$ which are not necessarily zero,
excluding $a_{\beta_1}$, which is taken to be $1$. So, now, we are
reduced to classify, modulo $\mathcal A$-equivalence, the
parametrizations in the set $\Sigma_\Gamma$, in order to classify
analytically plane branches.

Zariski noticed in \cite{[Z]} that, in order to get more
Elimination Criteria then the above ones, it was necessary to
consider the module of K\"ahler differentials over the local ring
of the branch, introducing new numerical analytic invariants.

Let us denote by
$$\Omega_2=\{ hdX+gdY; \ \ g,h\in \loc_2\},
$$
the $\loc_2$-free module of germs of differentials at $(\C^2,0)$,
and by
$$\Omega_1=\{ \xi dt; \ \ \xi\in \loc_1\},
$$
the $\loc_1$-free module of germs of differentials at $(\C,0)$.

A parametrization $\varphi:(\C,0) \to (\C^2,0)$, $t \mapsto
(x(t),y(t))$, induces a natural $\loc_2$-modules homomorphism
(thought as an extension of the map $\varphi^*:\loc_2 \to
\loc_1$):
\[
\begin{array}{rccc}
\varphi^*:& \Omega_2 & \longrightarrow & \Omega_1 .\\
       & hdX+gdY & \mapsto & (\varphi^*(h)x'(t)+\varphi^*(g)y'(t))dt
       \end{array}
       \]

The image of $\Omega_2$ in $\Omega_1$ under $\varphi^*$ is
isomorphic to the quotient of the module of K\"ahler differentials
over the local ring ($\varphi^*(\loc_2)$) of the branch determined
by $\varphi$, by its torsion submodule (cf. \cite{[HH1]} or
\cite{[HH3]}).

A primitive parametrization $\varphi$ also induces a valuation
$v_\varphi$ on $\Omega_2$, defined by
\[
v_\varphi(hdX+gdY)=\ord_t(\varphi^*(h)x'(t)+\varphi^*(g)y'(t))+1.
\]

We now define
\[
\Lambda_\varphi=v_\varphi(\Omega_2).
\]

The set $\Lambda_\varphi$ is $\mathcal A$-invariant, as we will
show later, and will play a key role in our solution of the
classification problem.

Since for all $h\in \M_2$ we have that
$v_\varphi(dh)=v_\varphi(h)$, it follows that
$\Gamma_\varphi\setminus\{0\} \subset \Lambda_\varphi$; and
because $\Gamma_\varphi$ has a conductor, we have that the set
$\Lambda_\varphi\setminus \Gamma_\varphi$, as a subset of the set
of gaps of $\Gamma_\varphi$, is finite.

Zariski in \cite{[Z1]} has shown that \ $\Lambda_\varphi \setminus
\Gamma_\varphi =\emptyset$ \ if, and only if, $\varphi$ is
$\mathcal A$-equivalent to the parametrization
$(t^{v_0},t^{v_1})$, where $\GCD(v_0,v_1)=1$.

Since $\Lambda_\varphi$ and $\Gamma_\varphi$ are $\mathcal
A$-invariant, it follows that, if \ $\Lambda_\varphi \setminus
\Gamma_\varphi \neq \emptyset$, \ then
$$
\lambda=\min\left( \Lambda_\varphi \setminus
\Gamma_\varphi\right)-v_0,$$  is an invariant under the $\mathcal
A$-action, called the Zariski invariant of $\varphi$. It is known
(see \cite{[Z]}) that such a $\varphi$ is $\mathcal A$-equivalent
to a Puiseux parametrization of the form
\begin{equation}
\varphi=(t^{v_0}, t^{v_1}+t^\lambda + \sum_{\lambda<i<c}a_it^i),
\end{equation}
and, in this case,
\[\lambda=v_\varphi(v_0XdY-v_1YdX)-v_0.
\]

 Related to the invariant $\lambda$, Zariski in \cite{[Z]}
proved the following extra elimination criterion:

\noindent {\bf EC3)} If $\varphi$ is as in (2.4) and $j-\lambda$
is in the semigroup generated by $v_0$ and $v_1$, then $\varphi$
is ${\mathcal A}$-equivalent to a para\-me\-tri\-zation $
(t^{v_0},t^{v_1}+t^\lambda + \sum_{\lambda<i<c}{a'_it^i}),$ with
$a'_i=a_i$, when $i<j$, and $a'_j=0$.

The above criterion doesn't work for all the equisingularity
class, but depends upon the $\mathcal A$-equivalence class of the
parametrization $\varphi$.

In the next theorem, our central result in this work, we will
determine all possible such elimination criteria, which will lead
us to what we call the {\em normal forms} for the Puiseux
parametrizations.

\begin{theo}[The Normal Forms Theorem]
Let $\varphi\in \Sigma_\Gamma$ be a Puiseux parametrization of a
plane branch with semigroup of values $\Gamma=\langle
v_0,v_1,\ldots,v_g\rangle$. Then, either $\varphi$ is $\mathcal
A$-equivalent to the monomial parametrization $(t^{v_0},t^{v_1})$,
or it is $\mathcal A$-equivalent to a parametrization
\begin{equation}
(t^{v_0},t^{v_1}+t^\lambda+\sum_{\stackrel{i>\lambda}{i\not \in
\Lambda-v_0}} a_it^i),
\end{equation}
where $\lambda$ is its Zariski invariant and
$\Lambda=\Lambda_\varphi$ is the set of orders of differentials of
the branch. Moreover, if $\varphi$ and $\varphi'$ {\rm (}with
coefficients $a_i'$ instead of $a_i${\rm )} are parametrizations
as in (2.5), representing two plane branches with same semigroup
of values and same set of values of differentials, then $\varphi
\sim_{\mathcal A} \varphi'$ if and only if there is $r\in \mathbb
C^*$ such that $r^{\lambda-v_1}=1$ and $a_i=r^{i-v_1}a_i'$, for
all $i$.
\end{theo}

Remark that the $\mathcal A$-normal form in (2.5) is completely
determined by the semigroup $\Gamma$ and the set $\Lambda$. So,
once $\Gamma$ is fixed, the number of ${\mathcal A}$-normal forms
is finite, corresponding to all possible sets $\Lambda$ in the
equisingularity class determined by $\Gamma$, which may be
computed by the algorithm presented in \cite{[HH3]}.

The above theorem gives the ultimate elimination criterion {\bf
EC} that contains all the known criteria {\bf EC1}, {\bf EC2} and
{\bf EC3}: \vspace{1mm}

\noindent {\bf EC)} If $\varphi$ is as in (2.4) and $j+v_0\in
\Lambda_\varphi$, $j> \lambda$, then $\varphi$ is ${\mathcal
A}$-equivalent to a parametrization $ (t^{v_0},t^{v_1}+ t^\lambda+
\sum_{\lambda<i<c}{a'_it^i}),$ with $a'_i=a_i$, when $i<j$, and
$a'_j=0$. \vspace{1mm}

The rest of the paper is devoted to prove Theorem 2.1.

\section{Orbits and their Tangent Spaces}

We will assume the reader familiar with the language of
singularity theory. We will denote by $j^k(h)$ the $k$-jet of an
element $h$.

Let Aut($\mathbb C^n,0$) denote the group of germs of analytic
automorphisms of $(\mathbb C^n,0)$, and let ${\rm Aut}_1(\mathbb
C^n,0)$ the subgroup of elements $A\in {\rm Aut}(\mathbb C^n,0)$
such that $j^1(A)={\rm Id}$.

We also denote by $\widetilde{{\rm Aut}}(\mathbb C^2,0)$ the
subgroup of elements $A$ of ${\rm Aut}(\mathbb C^2,0)$ such that
$j^1(A)=(X+\beta Y,Y)$, with $\beta\in \mathbb C$.

We say that the Puiseux parametrizations $\varphi_1$ and
$\varphi_2$ are ${\mathcal A}_1$-equivalent (resp.
$\widetilde{{\mathcal A}}$-equivalent) if $\varphi_2=\sigma\circ
\varphi_1 \circ \rho^{-1}$ with $\sigma\in {\rm Aut}_1(\mathbb
C^2,0)$ (resp. $\sigma\in \widetilde{{\rm Aut}}(\mathbb C^2,0)$)
and $\rho \in {\rm Aut}_1(\mathbb C,0)$.

We say that $\varphi_1$ and $\varphi_2$ in $\Sigma_\Gamma$ are
{\em homothetic}, or $\mathcal H$-equivalent if $\varphi_2=\sigma
\circ \varphi_1 \circ \rho^{-1}$, with $\rho(t)=\alpha t$ and
$\sigma(X,Y)=(\alpha^{v_0}X,\alpha^{v_1}Y)$, for some $\alpha\in
\mathbb C^*$.

So, the ${\mathcal A}$-action on the space of Puiseux series
representing an equisingularity class may be obtained by the
$\widetilde{{\mathcal A}}$-action followed by the $\mathcal
H$-action.

If ${\mathcal G}$ represents one of the actions ${\mathcal A}$,
${\mathcal A}_1$ or $\widetilde{{\mathcal A}}$, then ${\mathcal
G}^k$ will represent the Lie group action of $k$-jets of
corresponding automorphisms on the space $\Sigma_\Gamma^k$ of
$k$-jets of elements of $\Sigma_\Gamma$.

 Let us recall (a special
case) of the Complete Transversal Theorem of {\cite{[BKP]},
adapted to our use: \vspace {1mm}

\noindent {\bf The Complete Transversal Theorem.} \ {\em Let $G$
be a Lie group acting smoothly on an open set $U$ of an affine
space $A$ with underlying vector space $V$, and let $W$ be a
subspace of $V$ such that \ $\forall g\in G$, \ $\forall v\in U$
and $\forall w\in W$ with $v+w\in U$ and $g\cdot v+w\in U$, one
has
\[
g\cdot (v+w)=g\cdot v+w.
\]
If $W \subset T_v(G\cdot v)$, with $v\in U$, and $T_v(G\cdot v)$
is the tangent space at $v$ to the orbit $G\cdot v$, then for
every $w\in W$ such that $v+w\in U$, one has
$$ G(v+w)=G\cdot v.$$}

Although we are mainly interested in the ${\mathcal
A}$-equivalence, we will start analyzing the unipotent ${\mathcal
A}_1$-action, passing to the $\widetilde{{\mathcal A}}$-action
and, finally, applying homotheties, to get to the ${\mathcal
A}$-equivalence.

Let $U$ be the open set $\Sigma_\Gamma^k$ of the appropriate
affine space and let $G={\mathcal A}_1^k$. Then the initial
hypothesis at the beginning of the Complete Transversal Theorem is
fulfilled for $W=\{ (0,bt^k); b\in \mathbb C\}$.

The tangent spaces to the orbits ${\mathcal A}_1^k\cdot \varphi$
and $\widetilde{\mathcal A}^k\cdot \varphi$ at an element
$\varphi=(x(t),y(t)) \in \Sigma_\Gamma^k$ are given by:
\begin{equation}
 T_\varphi({\mathcal A}_1^k \cdot \varphi)=\left\{
j^k((x'(t),y'(t))\epsilon+(\varphi^*(g), \varphi^*(h) ); \
\epsilon \in \M^2_1, g,h\in \M^2_2\right\},
\end{equation}
and
\begin{equation}    T_\varphi(\widetilde{\mathcal A}^k \cdot
\varphi)= \left\{ j^k((x'(t),y'(t))\epsilon+(\varphi^*(g),
\varphi^*(h) ); \ \epsilon \in \M_1^2, h\in \M^2_2 ,\right.
\end{equation}
\[
\hspace{9cm} g\in \langle X^2,Y\rangle \left. \right\}.
\]

The proof of (3.1) may be found, for example, in \cite{[G]}, while
the proof of (3.2) may be obtained in a similar way.

 We will show how, by using the
Complete Transversal Theorem, one obtains all normal forms of
Puiseux parametrizations, with respect to the ${\mathcal
A}_1$-equivalence, by eliminating terms in the expansion of
$y(t)$, finding more elimination criteria than the general ones of
Ebey and Zariski, adapted to a specific branch. The idea is to
verify at each step if the $k$-jet of the parametrization is
${\mathcal A}_1^k$-equivalent to its $(k-1)$-jet, which implies
that the term of degree $k$ in $y(t)$ can be eliminated under the
${\mathcal A}_1$-action. For this, according to the Complete
Transversal Theorem, it is enough to verify if the vector
$(0,bt^k)$ belongs to the tangent space to the $\mathcal
A_1^k$-orbit of the $k$-jet of the parametrization, and this fact
may be expressed in terms of the existence of differentials in
$(\C^2,0)$ of certain order with respect to the valuation
determined by the parametrization, as we will see soon. The
procedure will stop after finitely many steps since all terms in
$y(t)$ of order greater or equal to the conductor $c$ of the
semigroup of values of the branch are eliminable. Next, we will
find the normal forms under the $\widetilde{{\mathcal A}}$-action
by analyzing separately some few remaining cases. Finally, the
normal forms under the ${\mathcal A}$-action are obtained applying
homotheties.

In order to apply this procedure, we will need to describe more
explicitly the tangent spaces to orbits in $\Sigma_\Gamma^k$.

\begin{lem} Let $k>v_1$ and $\varphi\in \Sigma_\Gamma^k$.
For $b\neq 0$, we have that the vector $(0,bt^k)$ belongs to
$T_{\varphi}({\mathcal A}_1^k \cdot \varphi)$ {\rm (}resp. to
$T_{\varphi}(\widetilde{\mathcal A}^k \cdot \varphi)${\rm )}, if
and only if there exist $g,h\in \M^2_2$ {\rm (}resp. $g\in \langle
X^2,Y\rangle$, $h\in \M^2_2${\rm )} such that
\begin{equation} \label{ordem}
k+v_0-1= {\rm ord}_t(\varphi^*(h)x'(t)- \varphi^*(g)y'(t)).
\end{equation}
\end{lem}
\noindent \pf We prove the result for $T_{\varphi}({\mathcal
A}_1^k \cdot \varphi)$, since the other situation is similar.

In order to have $(0,bt^k)\in T_{\varphi}(\mathcal A_1^k\cdot
\varphi)$ it is necessary and sufficient to be able to solve the
system:
$$\left \{\begin{array}{l}
x'(t).\epsilon+ \varphi^*(g) =0 \ \ \ \ \ mod\ t^{k+1}\\
y'(t).\epsilon+ \varphi^*(h)=bt^k \ \ \ mod\
t^{k+1}.\end{array}\right .$$

That is, $\displaystyle \epsilon =-\frac{\varphi^*(g)}{x'(t)}\
mod\ t^{k+1}$. Notice that $\epsilon\in\M^2_1$, since
$g\in\M^2_2$.

In this way we get the equation
$$bt^k=\frac{\varphi^*(h)x'(t)-\varphi^*(g)y'(t)}{x'(t)}\ \  mod\ t^{k+1}.$$

So, $(0,bt^k)\in T_{\varphi}(\mathcal A_1^k\cdot \varphi)$ if, and
only if, there exist $g,h\in\M^2_2$ satisfying (\ref{ordem}). \fim

For $i\in \N$, we define
$$\Omega_2^{(i)}= \{hdX+gdY\in \Omega_2; \ \ g, h \in \M_2^i \},
$$
where we put $\M_2^0=\loc_2$. So, $\Omega_2^{(0)}=\Omega_2$.

Given a primitive parametrization $\varphi$, we also define
$$\Lambda_\varphi^{i}=v_\varphi(\Omega_2^{(i)}).$$
Notice that $\Lambda_\varphi^{0}=\Lambda_\varphi$. These sets are
invariant under $\mathcal A$-equivalence, as we show below.

\begin{prop} If $\varphi$ and $\varphi_1$ are $\mathcal A$-equivalent
primitive parametrizations, then, for all $i\in \N$, we have
$\Lambda_\varphi^{i}=\Lambda_{\varphi_1}^{i}$.
\end{prop}
\noindent \pf The commutative diagram,
\[
\begin{array}{ccc}
   \C,0 & \stackrel{\varphi}{\longrightarrow} & \C^2,0 \\
   \downarrow \rho & &                  \downarrow
   \sigma\\
\C,0 & \stackrel{\varphi_1}{\longrightarrow} & \C^2,0
\end{array}
\]
with isomorphisms $\sigma$ an $\rho$ induces the following two
diagrams:
\[
\begin{array}{lclclcl}
   \loc_2 & \stackrel{\varphi^*}{\longrightarrow} & \loc_1 & \quad &
   \Omega_2 & \stackrel{\varphi^*}{\longrightarrow} & \Omega_1\\
   \uparrow \sigma^* &  &                  \uparrow
   \rho^* & \quad & \uparrow \sigma^* &  &                  \uparrow
   \rho^*  \\
\loc_2 & \stackrel{\varphi_1^*}{\longrightarrow} & \loc_1 & \quad
&  \Omega_2 & \stackrel{\varphi_1^*}{\longrightarrow} & \Omega_1
\end{array}
\]
with isomorphisms $\sigma^*$ and $\rho^*$. Now the result follows
by functoriality, observing that $\sigma^*(\M_2^{i})=\M_2^{i}$.
\fim

 If we define $$\Omega_2'=\{hdX+gdY\in \Omega_2; \
\ g\in \langle X^2,Y \rangle, h\in \M_2^2\},$$ and
$\Lambda_\varphi'=v_\varphi(\Omega_2')$, then Lemma 3.1 may be
rephrased as follows:
\begin{prop}
 Let $k>v_1$ and $\varphi\in \Sigma_\Gamma^{k}$.
For $b\neq 0$, we have that $(0,bt^k)$ belongs to
$T_{\varphi}({\mathcal A}_1^k \cdot \varphi)$ {\rm (}resp. to
$T_{\varphi}(\widetilde{\mathcal A}^k \cdot \varphi)${\rm )} if,
and only if,
\begin{equation} \label{ordemdif}
k+v_0\in \Lambda_\varphi^{2} \ \ \hbox{{\rm (}resp. $k+v_0\in
\Lambda_\varphi'${\rm )}}
\end{equation}
\end{prop}

\section{Normal $\mathcal A_1$-Forms}

In this section we will find the normal forms of Puiseux
parametrizations in an equisingularity class with given semigroup
of values $\Gamma$, under $\mathcal A_1$-equivalence. We begin
with a proposition that will give us the recursion step.
\begin{prop}
Let $\varphi=( t^{v_0}, t^{v_1} + \sum_{v_1<i<c} a_it^i )\in
\Sigma_\Gamma$ and let $k$ be an integer such that $k+v_0\in
\Lambda_\varphi^{2}$. Then there exists $\varphi_1\in
\Sigma_\Gamma$ such that $\varphi_1 \sim_{\mathcal A_1} \varphi$
and
\[ j^k(\varphi_1)=j^{k-1}(\varphi_1)=j^{k-1}(\varphi).\]
\end{prop}
\noindent \pf From Proposition 3.3, we have that the vector
$(0,-a_kt^k)$ belongs to $T_{j^{k}(\varphi)}({\mathcal A}_1^k\cdot
j^{k}(\varphi))$, and therefore by the Complete Transversal
Theorem it follows that $j^k(\varphi) \sim_{\mathcal A_1^k}
j^{k-1}(\varphi)$. Hence, there are appropriate germs of analytic
isomorphisms $\sigma$ and $\rho$ such that $\sigma \circ
j^k(\varphi) \circ \rho^{-1}=j^{k-1}(\varphi)$. So, $j^k(\sigma
\circ \varphi \circ \rho^{-1})=j^{k-1}(\varphi)$. Now the result
follows putting $\varphi_1=\sigma \circ \varphi \circ \rho^{-1}$
\fim

The following result will be important in the sequel.
\begin{prop} Let $\varphi=(t^{v_0},t^{v_1}+t^{\lambda}+\cdots)$ be a Puiseux
parametrization with $\Gamma_\varphi=\langle
v_0,v_1,\ldots,v_g\rangle$. If
$S=\{v_0,2v_0,v_1,v_0+v_1,2v_1,v_0+\lambda\}$, then one has
\[
S \ \ \subseteq \ \ \Lambda_\varphi \setminus \Lambda_\varphi^{2}
\ \ \subseteq \ \ S  \cup  \{v_1+\lambda\},
\]
with equality on the top if, and only if, $n_1=2$ and $g\geq 2$.
\end{prop}
\noindent \pf We have that $n\in \Lambda_\varphi\setminus
\Lambda_\varphi^2$ if, and only if, $n=v_\varphi(\omega)$, where
 $\omega=hdX+gdY$, with $h\not\in \M_2^{2}$ or $g\not\in \M_2^{2}$.

We have that $S \ \subseteq \Lambda_\varphi \setminus
\Lambda_\varphi^{2}$, \ since $v_\varphi(dX)=v_0$,
$v_\varphi(dY)=v_1$, $v_\varphi(XdX)=2v_0$,
$v_\varphi(XdY)=v_0+v_1$, $v_\varphi(YdY)=2v_1$, and
$v_\varphi(v_1YdX-v_0XdY)=v_0+\lambda$.

Now, suppose that $v_\varphi(hdX+gdY)\not \in S$, where $h=\alpha
X+\beta Y+h_2$ and $g=a X+ b Y+g_2$, with $h_2,g_2\in \M_2^2$, and
one of the numbers $a,b, \alpha$ or $\beta$ is not zero.

So, in this case, we must have
$v_\varphi(hdX+gdY)>v_\varphi(hdX)=v_\varphi(gdY)$. This implies
that
\[
v_\varphi(h)+v_0=v_\varphi(g)+v_1,
\]
with $v_\varphi(h)=v_0$ or $v_\varphi(h)=v_1$ or
$v_\varphi(g)=v_0$ or $v_\varphi(g)=v_1$.

If $v_\varphi(h)= v_0$, then we would have
$2v_0=v_\varphi(g)+v_1$, which is not possible. Hence $\alpha=0$.

If $v_\varphi(g)=v_0$, then $v_\varphi(h)+v_0=v_1+v_0$, hence
$v_\varphi(h)=v_1$.

If $v_\varphi(h)=v_1$, then $v_1+v_0=v_\varphi(g)+v_1$, so
$v_\varphi(g)=v_0$. Hence, in this case,
$v_\varphi(hdX+gdY)=v_0+\lambda \in S$, which is to be excluded.
Hence, $a=\beta=0$.

So, the only remaining possibility is that $v_\varphi(g)=v_1$, in
which case, $b\neq 0$ and $a=\alpha=\beta= 0$. So, we have
$v_\varphi(h)+v_0=2v_1$, hence $v_1<v_\varphi(h)<2v_1$, which
implies that $v_\varphi(h)=sv_1+rv_0$, with $s=0,1$. We have that
$s=0$, because, otherwise, we would have $v_1+(r+1)v_0=2v_1$,
which would imply that $v_0$ divides $v_1$, a contradiction.

If the genus of $\varphi$ is 1, then $v_0=2$, and in this case,
because of {\bf EC1} we have that $\varphi$ is $\mathcal
A$-equivalent to the parametrization $(t^2,t^{v_1})$, hence not
satisfying the hypothesis of the Proposition.

Therefore, $g\geq 2$, and $v_\varphi(h)=rv_0$. Therefore,
\[
(r+1)v_0=2v_1,
\]
which implies $n_1=2$. Also, $v_\varphi(hdX+gdY)>2v_1$, which in
view of the expression of $\varphi$ and the above equality implies
that $v_\varphi(hdX+gdY)=v_1+\lambda$.

Conversely, if $g\geq 2$ and $n_1=2$, we have that
\[
v_\varphi(v_1X^rdX-v_0YdY)=v_1+\lambda,
\]
where $(r+1)v_0=2v_1$. \fim

Now, we have the following result:

\begin{prop} Let $\varphi\in \Sigma_\Gamma$ and set
$\Lambda=\Lambda_\varphi$. Suppose that $\Lambda\setminus
\Gamma\neq \emptyset$, and let $\lambda$ be the Zariski invariant
of $\varphi$. Then $\varphi$ is ${\mathcal A}_1$-equivalent to a
parametrization
\[
(t^{v_0},t^{v_1}+t^\lambda+ \sum_{\stackrel{i\not \in
\Lambda^{2}-v_0}{i>\lambda}}a_it^i).
\]
\end{prop}
\noindent \pf First observe that since $\Lambda\setminus
\Gamma\neq \emptyset$, it follows that $v_0\geq 3$, so, in this
situation, any integer $l+v_0$, where $l$ is greater or equal than
the conductor $c$ of $\Gamma$, belongs to $\Lambda$ but, by
Proposition 4.2, it is not in $\Lambda \setminus \Lambda^{2}$, so
it is in $\Lambda^{2}$. This shows that the set $\N \setminus
(\Lambda^{2}-v_0)$ is finite (bounded by above by $c-1$).

Let $\lambda_1,\ldots, \lambda_s$ be the elements in
$\Lambda^{2}-v_0$ in the interval $(\lambda, c)$. From Proposition
4.1, there exists a Puiseux parametrization $\varphi_1$ with
$\varphi_1 \sim_{\mathcal A_1} \varphi$ such that
\[
j^{\lambda_1}(\varphi_1)=j^{\lambda_1-1}(\varphi_1)= \
j^{\lambda_1-1}(\varphi).
\]
Next, do the same with $\varphi_1$ instead of $\varphi$ and
$\lambda_2$ instead of $\lambda_1$, observing, by Proposition 3.2,
that $\Lambda_{\varphi_1}^{2}=\Lambda^{2}$; etc. \fim

The next step will be to pass from the ${\mathcal
A}_1$-equivalence to the $\widetilde{\mathcal A}$-equivalence.

\section{Passage from the ${\mathcal
A}_1$-equivalence to the $\widetilde{\mathcal A}$-equivalence}

To get the normal forms of Theorem 2.1, in view of the result of
Proposition 4.3, it suffices to show that the terms in $y(t)$ of a
Puiseux parametrization of order $k$ such that $k>\lambda$ and
$k\in (\Lambda\setminus \Lambda^{2})-v_0$, may be eliminated
without changing the preceding terms.

Terms of order $k\in S-v_0$, where $S$ is as in Proposition 4.2,
excepting $k=\lambda$, may be eliminated by {\bf EC2}. The only
remaining possibility are terms of order $v_1+\lambda-v_0$, when
$g\geq 2$ and $n_1=2$ (cf. Proposition 4.2), which we will show
below how to eliminate them without changing the preceding ones.

To do this, we will need to analyze more closely the
$\widetilde{\mathcal A}$-action on Puiseux parametrizations.

Let $\varphi(t)=(t^{v_0},y(t))$, where $y(t)=\sum_{i}a_it^i$, and
let $\sigma$ and $\rho$ as in (2.2), but with $r=1$, $p=\beta
Y+p_1$, where $\beta\in \C$ and $p_1,q\in \M_2^2$.

Now, considering the expression of $\rho$ in (2.2), raising both
sides to the power $i$ and then applying the binomial expansion we
get
\[
t_1^i=t^i\left[ \sum_{j=0}^\infty {i/v_0 \choose j} \left(
\frac{p(t)}{t^{v_0}}\right) ^j\right],
\]
where $p(t)=\varphi^*(p)$.

By using this in the expression $y(t)=\sum_ia_it^i$, we get
\[
y(t_1)=y(t)+\sum_{i}a_it^i\frac{i}{v_0}\frac{p(t)}{t^{v_0}} +
A(t),
\]
where
\[
A(t)=\sum_ia_it^i\sum_{j=2}^\infty {i/v_0 \choose
j}\left(\frac{p(t)}{t^{v_0}} \right)^j.
\]

Now, from the expression of $y_1(t_1)$ in (2.3) we get
\begin{equation}
 y_1(t_1)=y(t_1)+B(t),
\end{equation}
where, if we put $q(t)=\varphi^*(q)$,
\begin{equation}
B(t)=\frac{q(t)x'(t)-p(t)y'(t)}{x'(t)}-A(t).
\end{equation}

\begin{prop} Let $\varphi(t_1)=(t_1^{v_0}, y(t_1))$, where
$y(t_1)=t_1^{v_1}+t_1^{\lambda}+\sum_{i>\lambda}a_it_1^i$, be a
Puiseux parametrization, such that the genus of $\varphi$ is
greater than $1$, and $n_1=2$. Then there exists
$y_1(t_1)=t_1^{v_1}+t_1^{\lambda}+\sum_{i>\lambda}a'_it_1^i$, with
$a'_i=a_i$, for $i<v_1+\lambda-v_0$ and $a'_{v_1+\lambda-v_0}=0$,
such that $(t_1^{v_0},y_1(t_1))\sim_{\widetilde{\mathcal A}}
(t_1^{v_0},y(t_1))$.
\end{prop}
\noindent \pf Let $\beta\in \C$ and $p_1,q\in \M_2^2$, as above.
We will show that we may choose $p_1$ and $q$, in such a way that
$\ord_{t_1}(B(t))=v_1+\lambda-v_0$, where $B$ is as in (5.2), and
then, by adjusting the value of $\beta$, we may make this term
cancel the corresponding one in $y(t_1)$ in equation (5.1).

Since $t_1=\rho(t)$, with $\rho$ an automorphism of $\loc_1$, we
have that
\[
\ord_{t_1}(B(t))=\ord_{t}(B(t)),
\]
so, we may work with the powers of $t$ in the expression of
$B(t)$.

Remark that $n_1=2$ implies that $m_1v_0=2v_1$, where
$m_1=\beta_1/e_1=v_1/e_1$.

Now we choose $p_1=0$ and $q=\frac{v_1}{v_0}\beta X^{m_1-1}+g$,
with $g\in \M_2^2$ such that $v_\varphi(g)>(m_1-1)v_0$.

Let us write
\[
B(t)=B_0(t)+B_1(t)+B_2(t),
\]
where
\[
B_0(t)= \beta \frac{(v_1/v_0) x(t)^{m_1-1}x'(t)-y(t)y'(t)}{x'(t)},
\]
\[
B_1(t)=\frac{g(t)x'(t)}{x'(t)}=g(t),
\]
and
\[
B_2(t)=-A(t)=-\sum_{i\geq v_1}a_it^i\left( \sum_{j=2}^{\infty}
{i/v_0 \choose j}\left( \frac{\beta
y(t)}{t^{v_0}}\right)^j\right).
\]

A direct computation shows that, if $\beta\neq 0$, then
\[
v_\varphi(B_0)=v_\varphi((v_1/v_0)
x(t)^{m_1-1}x'(t)-y(t)y'(t))-(v_0-1)=
\]
\[\
v_1+\lambda-1-v_0+1=v_1+\lambda-v_0.
\]

On the other hand, by expanding $B_2(t)$ one sees that there will
be terms either of degree greater than $v_1+\lambda-v_0$, or of
degree $rv_0+sv_1$, greater than $(m_1-1)v_0$, which can be
eliminated by a suitable choice of $g(t)$. \fim

With this last proposition we finished the proof of the existence
part of Theorem 2.1, concerning the normal forms.

Now, to prove that if two Puiseux parametrizations are $\mathcal
A$-equivalent, then they are conjugate under homothety, it will be
sufficient to prove that if two Puiseux parametrizations are
$\widetilde{\mathcal A}$-equivalent, then they are equal, because
the $\mathcal A$-action is decomposed into the
$\widetilde{\mathcal A}$-action and the $\mathcal H$-action.

Fixing a set $\Lambda$ of values of differentials in the
equisingularity class determined by a semigroup $\Gamma$, let us
consider the linear space
\[
N_\Lambda =\{ (t^{v_0},t^{v_1}+t^{\lambda}+\sum_{j>\lambda}a_jt^j
)\in \Sigma_\Gamma; \ a_j=0, \ \hbox{for $j\in \Lambda-v_0$} \}.
\]

If we denote by $N^k_\Lambda$ the space $j^k(N_\Lambda)$, we have
the following lemma:

\begin{lem} If $\alpha\in N_\Lambda$, then for all $k>\lambda$, we
have
\[
N^k_\Lambda \cap T_{j^k(\alpha)}(\widetilde{\mathcal A}^k\cdot
j^k(\alpha))=\{ j^k(\alpha)\}.
\]
\end{lem}
\noindent \pf  Suppose the assertion not true. Take $k$ minimal
with the following property:
\[
N^k_\Lambda \cap T_{j^k(\alpha)}(\widetilde{\mathcal A}^k\cdot
j^k(\alpha))\neq\{ j^k(\alpha)\}.
\]

So, there exists $\beta\in N^k_\Lambda \cap
T_{j^k(\alpha)}(\widetilde{\mathcal A}^k\cdot j^k(\alpha))$ such
that $\beta\neq j^k(\alpha)$ and $j^{k-1}(\beta)=j^{k-1}(\alpha)$.
Therefore, there exists $b\in \C^*$ such that
\[
\beta-j^k(\alpha)=(0,bt^k)\in T_{j^k(\alpha)}(\widetilde{\mathcal
A}^k\cdot j^k(\alpha)).
\]

Hence, from Proposition 3.3, it follows that $k\in \Lambda-v_0$.
But, since $j^k(\alpha)\in N_\Lambda^k$, it follows that
$j^k(\alpha)=j^{k-1}(\alpha)$. So, for some $b\neq 0$,
\[
\beta=j^{k-1}(\alpha)+(0,bt^k).
\]

But, since $\beta\in N_\Lambda^k$, one should have $b=0$, which is
a contradiction. \fim

Now we proceed to prove the uniqueness of the $\widetilde{\mathcal
A}$-normal forms.

Let
$\varphi(t)=(t^{v_0},t^{v_1}+t^\lambda+\sum_{j>\lambda}a_jt^j)\in
\Sigma_\Gamma$ be a Puiseux parametrization with
$\Lambda_\varphi=\Lambda$. We denote by $\widetilde{\mathcal
A}^{c-1}\cdot \varphi$ the orbit of $\varphi$ in $\Sigma_\Gamma$,
with respect to the $\widetilde{\mathcal A}^{c-1}$-action.

We want to show that
\[
N_\Lambda\cap \widetilde{\mathcal A}^{c-1}\cdot
\varphi=\{\varphi\}.
\]

Indeed, if $N_\Lambda\cap \widetilde{\mathcal A}^{c-1}\cdot
\varphi\neq\{\varphi\}$, take $\varphi_1\in N_\Lambda\cap
\widetilde{\mathcal A}^{c-1}\cdot \varphi$, with $\varphi_1\neq
\varphi$. Since $\widetilde{\mathcal A}^{c-1}\cdot \varphi$ is
arcwise connected, there exists an arc in $\widetilde{\mathcal
A}^{c-1}\cdot \varphi$ joining $\varphi$ to $\varphi_1$. Since
reduction to the normal form is continuous, it follows that
$\varphi$ wouldn't be an isolated point in $N_\Lambda
\cap\widetilde{\mathcal A}^{c-1}\cdot \varphi$. But this is a
contradiction because of Lemma 5.2.

\section{Zariski's Problem and Computational Aspects}

Our methods, more than describing all normal forms (up to
homotheties) of plane branches, with respect to analytic
equivalence, give an effective way to obtain the normal form of a
given branch and, as well, to distinguish from analytic point of
view two given plane branches.

Indeed, since the set $\Lambda$ is an analytic invariant, it is
the same for the branch and its normal form, so given two
parametrizations $\varphi_1$ and $\varphi_2$ with same semigroup
of values, we compute with the procedures exposed in \cite{[HH3]}
(specially Algorithm 4.10) the sets $\Lambda$ for both. If these
are distinct, the branches are not equivalent. If they are equal,
we proceed to put the parametrizations under they normal forms. To
put a given parametrization $\varphi$ into its normal form, it is
enough to consider changes of coordinates corresponding to an
element of the group $\widetilde{\mathcal A}^{l-v_0}$, where $l$
is the greatest integer not in $\Lambda$. More precisely, taking
\[
p=\sum_i \alpha_i\prod_{j=0}^{g}h_j^{a_{ij}}, \ \ \ q=\sum_i
\gamma_i\prod_{j=0}^{g}h_j^{b_{ij}},
\]
where $h_0,h_1,\ldots,h_g$ are elements of a standard basis for
the local ring ${\mathcal O}_f$ of the branch $(f)$ associated to
$\varphi$ (see \cite{[HH3]} for the definition (Definition 2.1)
and how to compute them (Algorithm 3.2)), and the $\alpha_i$'s and
$\gamma_i$'s are parameters, such that in the development in power
series the smallest order term of $p$ is greater than $v_0$ and of
$q$ is greater than $v_1$ and the terms of orders belonging to
$\Gamma$ in $p$ (resp. $q$) are less than $l-v_1$ (resp. less than
$l-v_0$).

Performing an action as (2.2) with the $p$ and $q$ as above, which
is computationally possible, we impose conditions on the
coefficients $\alpha_i$ and $\gamma_i$ in order to bring the
parametrization into its normal form, in the way we did during the
proof of Proposition 5.1.

This done, the analytic equivalence reduces to verify homothety,
which is trivial.

Let us remark that Ebey (cf. \cite{[E]}, Theorem 5), by using
arguments from the theory of algebraic groups, predicted the
existence of some kind of normal forms under the $\mathcal
A_1$-equivalence (cf. our Proposition 4.3) which he called {\em
canonical forms}, but without any indication on how they could be
obtained nor how they should look like.

The {\em stratified moduli problem} is also solved, since it is
the disjoint union of quotients of a finite number of
semi-algebraic sets, by finite groups, corresponding to the
quotients modulo finitely many homotheties of the normal form
corresponding to a given $\Lambda$ under the $\widetilde{{\mathcal
A}}$-equivalence. The sets $\Lambda$ and the conditions on the
coefficients that determine them, for a fixed equisingularity
class, may be computed by the algorithms we developed in
\cite{[HH3]}. One of these set $\Lambda_{gen}$ corresponds to the
generic branch, easily recognized by the open conditions on the
coefficients.

Finally, the dimension of the component of the moduli
corresponding to a given set $\Lambda$ is determined by the normal
form and is at most equal to the number of gaps of $\Lambda$
greater that $\lambda$, since some of the coefficients of the
parametrization may be fixed constants. The dimension of the
generic component is exactly equal to the number of gaps of
$\Lambda_{gen}$ greater that $\lambda$, since in this case, no
coefficient in the corresponding normal form is constant.

\section{Some examples}

In what follows we give two concrete examples of the application
of our method. The first example will describe a result obtained
in \cite{[HH2]}, and the second one is a new example which we
relate to a question posed by Zariski in \cite{[Z]}.

\begin{example} The table below gives the analytic classification
of plane branches in the equisingularity class of $\Gamma=\langle
6,9,19 \rangle$.
\begin{center}
\begin{tabular}{|l|l|}
\hline Condition  & Normal Form \\
\hline $ b\not\in\{-\frac{1}{2},\frac{29}{18}\} $ & $(t^6,
t^9+t^{10}+bt^{11}+b_1t^{14}+b_2t^{17})$
\\ \hline $b=\frac{29}{18}$ & $(t^6,
t^9+t^{10}+bt^{11}+b_1t^{14}+b_2t^{17}+b_3t^{23}) $ \\ \hline
$\begin{array}{c} b=-\frac{1}{2} \\
A\not =0 \end{array}$ & $(t^6,
t^9+t^{10}+bt^{11}+b_1t^{14}+b_2t^{17}+b_3t^{20})$ \\ \hline
$\begin{array}{c} b=-\frac{1}{2}
\\ A=0
\end{array}$ & $(t^6,
t^9+t^{10}+bt^{11}+b_1t^{14}+b_2t^{17}+b_3t^{20}+b_4t^{26}) $ \\
\hline
\end{tabular}
\vspace{2mm}
\end{center}

\noindent Where $A=14+\frac{769}{2}b_1-532b_2-576b_1^2$.
Moreover, two branches belonging to the same normal form are
equivalent if, and only if, they are equal.
\end{example}

In the above example we have that the stratified moduli space
consists of four strata, all of them of dimension $3$, with the
first one corresponding to the generic stratum.

\begin{example} Our second example deals with the classification of the
equisingularity class given by the semigroup of values
$\Gamma=\langle 7,8 \rangle$.

Since the conductor of $\Gamma$ is $42$, we have  $$
\Sigma_\Gamma= \{ (t^7,t^8+\sum_{8<i<42} a_it^i); \ a_i\in \C\}.$$

Algorithm 4.10 of \cite{[HH3]} and Theorem 2.1 give the following
table:

{\footnotesize
\begin{center}
\begin{tabular}{|c|l|c|}\hline
 Condition & Normal Form & $\Lambda\setminus\Gamma$ \\ \hline
$a_{12}\neq \frac{13+9a_{11}^2}{8}$ & $ (t^7,
t^8+t^{10}+a_{11}t^{11}+a_{12}t^{12}+a_{13}t^{13}+a_{20}t^{20}) $
& $\begin{array}{c} 17,25,26\\ 33,34,41\end{array}$ \\ \hline
$a_{13}\neq \frac{39}{10}a_{11}+\frac{27}{20}a_{11}^3$ & $(t^7,
t^8+t^{10}+a_{11}t^{11}+\frac{13+9a_{11}^2}{8}t^{12}+a_{13}t^{13}+a_{19}t^{19})
$ & $\begin{array}{c} 17,25,27\\ 33,34,41 \end{array}$ \\ \hline
$a_{20}\neq B$ & $\begin{array}{l} (t^7,
t^8+t^{10}+a_{11}t^{11}+\frac{13+9a_{11}^2}{8}t^{12}+ \\
\hspace{2cm} \left (
\frac{39}{10}a_{11}+\frac{27}{20}a_{11}^3\right
)t^{13}+a_{19}t^{19}+a_{20}t^{20}) \end{array}$ &
$\begin{array}{c}
17,25,33\\ 34,41\end{array}$ \\
\hline $a_{20}= B$ & $
\begin{array}{l}
(t^7, t^8+t^{10}+a_{11}t^{11}+\frac{13+9a_{11}^2}{8}t^{12}+\\
\hspace{.5cm} +\left (
\frac{39}{10}a_{11}+\frac{27}{20}a_{11}^3\right
)t^{13}+a_{19}t^{19}+a_{20}t^{20}+a_{27}t^{27})
\end{array}$ & $\begin{array}{c} 17,25\\ 33,41\end{array}$ \\ \hline
 & $(t^7, t^8+t^{11}+a_{12}t^{12}+a_{13}t^{13}+a_{20}t^{20})$
 & $\begin{array}{c} 18,25,26\\ 33,34,41\end{array}$ \\ \hline
 & $(t^7,t^8+t^{12}+a_{13}t^{13}+a_{18}t^{18})
$ & $\begin{array}{c} 19,26,27\\ 33,34,41\end{array}$ \\ \hline
$a_{18}\neq -\frac{1}{2}$ &
$(t^7,t^8+t^{13}+a_{18}t^{18}+a_{19}t^{19}+a_{26}t^{26}) $ &
$\begin{array}{c} 20,27,33\\ 34,41\end{array}$ \\ \hline
 & $(t^7, t^8+t^{13}-\frac{1}{2}t^{18}+a_{19}t^{19}+a_{26}t^{26})
$ & $\begin{array}{c} 20,27\\ 34,41\end{array}$ \\ \hline
$a_{20}\neq \frac{121}{120}a_{19}^2$&
$(t^7,t^8+t^{18}+a_{19}t^{19}+a_{20}t^{20}+a_{27}t^{27}) $ &
$\begin{array}{c} 25,33\\ 34,41\end{array}$ \\ \hline & $(t^7,
t^8+t^{18}+a_{19}t^{19}+\frac{121}{120}a_{19}^2t^{20}+a_{27}t^{27})
$ & $\begin{array}{c} 25,33,41\end{array}$ \\ \hline
 & $(t^7, t^8+t^{19}+a_{20}t^{20})
$ & $\begin{array}{c} 26,33\\ 34,41\end{array}$ \\ \hline
 & $(t^7, t^8+t^{20}+a_{26}t^{26})
$ & $27,34,41$ \\ \hline
 & $(t^7, t^8+t^{26}+a_{27}t^{27})
$ & $33,41$ \\ \hline
 & $(t^7, t^8+t^{27})
$ & $34,41$ \\ \hline
 & $(t^7, t^8+t^{34})
$ & $41$ \\ \hline & $ (t^7, t^8) $ & $\emptyset$ \\ \hline
\end{tabular}
\end{center}
}

\noindent where
$$B=\frac{11}{4}a_{11}a_{19}-
\frac{357}{512}-\frac{47399}{2560}a_{11}^2-\frac{10097}{320}a_{11}^4-\frac{17523}{1280}a_{11}^6
-\frac{2187}{1280}a_{11}^8.$$

Two parametrizations of the above table on the same line are
equivalent if and only if they are homothetic, with respect to the
appropriate root of unity {\rm ((}$\lambda-8${\rm )}-th root of
unity{\rm )}.
\end{example}

From the above table we see that the generic component of the
moduli, corresponding to parametrizations on the first line, has
dimension $4$. There are six strata of dimension $3$, three strata
of dimension $2$, three strata of dimension $1$ and three strata
of dimension $0$.

Zariski dedicated Sections $4,5$ and $6$ of Chapter VI, in
\cite{[Z]}, to the study of branches with semigroups of the form
$\Gamma =\langle v_0,v_0+1\rangle$, where the following result is
proved:

\noindent {\sc Theorem} \ (\cite{[Z]}, Th\'eor\`eme 6.12) {\em Let
$v_0\geq 5$, and for all $s\in\{ 2,\ldots ,v_0-2\}$, define
\[ {\mathcal L}_s=\{ sv_0+s+2,sv_0+s+3,\ldots ,sv_0+s+v_0-1-s\}.
\]
Let
\[
\varphi=(t^{v_0}, t^{v_0+1}+a_{v_0+3}t^{v_0+3}+\cdots
+a_{2v_0-1}t^{2v_0-1}+\displaystyle{\sum_{i\in\cup_{s=2}^{q}{\mathcal
L}_s}{a_it^i}}), \] where $q=\left [ \frac{v_0-3}{2}\right ]$ and
$a_i=0$, whenever $i$ is one of the first $s+1$  elements of
${\mathcal L}_s$, for all $2\leq s\leq q$.

Then two generic parametrizations of the above form are $\mathcal
A$-equivalent if, and only if, they homothetic.} \vspace{2mm}

In \cite{[Z]}, Zariski remarks that the above theorem is true for
$2\leq v_0\leq 6$ without the condition on the genericity of the
parameters (Remarque 6.14), and asks the following question:
\vspace{2mm}

{\it Is the above theorem true without the assumption of the
genericity on the coefficients of the parametrizations? }
\vspace{0.2cm}

The answer is no! And an example may be given considering branches
with semigroup $\Gamma =\langle 7,8\rangle$.

Consider
$$\varphi = (t^7, t^8+t^{10}+t^{11}+\frac{11}{4}t^{12}+a_{13}t^{13}+a_{20}t^{20}),
$$ with $a_{13}\not =\frac{21}{4}$. Obviously, $\varphi$ is in the form of the above theorem.

If we consider changes of coordinates as in (2.2) with
$$p=b_1x^2-\frac{3}{2}b_1xy-\frac{1}{4}b_1y^2+b_2x^3+b_3x^2y+b_4xy^2+b_5y^3$$
and
$$q=\frac{8}{7}b_1xy-\frac{12}{7}b_1y^2+\left ( -\frac{135}{28}b_1-\frac{15}{14}a_{13}b_1\right )
x^3+b_6x^2y+b_7xy^2+b_8y^3,$$ where
$$\begin{array}{l}
b_2=\frac{4}{7}b_1^2-\frac{227}{6}b_1+\frac{199}{24}a_{13}b_1-\frac{2}{3}b_3-\frac{8}{3}b_5, \\
b_3=6b_4-\frac{45}{2}b_1a_{13}^2+\frac{9565}{16}b_1-40b_8+\frac{72}{7}b_1^2+\frac{4297}{16}a_{13}b_1, \\
b_4=(4a_{13}-21)^{-1}\frac{1}{1120}
(-2720a_{13}b_1^2+557690a_{13}b_1- \\
\hspace{0.75cm} 277200b_1a_{13}^2+

+16800b_1a_{13}^3+26880a_{13}b_5+\\
\hspace{0.75cm} 2459289b_1-141120b_5+14280b_1^2+2688b_1a_{20}), \\
b_6=-\frac{25}{4}b_1-\frac{29}{28}a_{13}b_1+\frac{8}{7}b_2+\frac{4}{49}b_1^2, \\
b_7=-\frac{3}{2}a_{13}b_1+\frac{8}{7}b_3-\frac{641}{56}b_1-\frac{12}{49}b_1^2, \\
b_8=-\frac{52}{7}b_1+\frac{8}{7}b_4-\frac{81}{28}a_{13}b_1+\frac{2}{7}b_2+\frac{18}{49}b_1^2, \\
\end{array}$$
we have
$$\sigma\circ \varphi \circ \rho^{-1}(t_1)= ( t_1^7, t_1^8+t_1^{10}+t_1^{11}+\frac{11}{4}t_1^{12}+a_{13}t_1^{13}+
( a_{20}+5b_1 ( \frac{3}{4}-\frac{1}{7}a_{13} ) )t_1^{20}).
$$

By choosing conveniently $b_1$ we see that the normal form of
$\varphi$ is given in the second row of the table in Example 7.2,
but $\varphi$ itself is not in normal form (this gives the
reduction of $\varphi$ to normal form).

Since $a_{13}\not =\frac{21}{4}$, then for each $b_1\not =0$ we
get a parametrization $\mathcal A$-equivalent to $\varphi$, as
described in Zariski's Theorem, without being homothetically
equivalent to $\varphi$, giving a negative answer to Zariski's
question.

\vspace{5mm}

\noindent Abramo Hefez\\
Universidade Federal Fluminense\\
Instituto de Matem\'atica \\
R. Mario Santos Braga, s/n \\
24020-140 Niter\'oi, RJ - Brazil\\
{\small {\em E-mail}:\quad hefez@mat.uff.br}

\vspace{3mm}

\noindent Marcelo E. Hernandes\\
Universidade Estadual de Maring\'a \\
Departamento de Matem\'atica\\
Av. Colombo, 5790  \\
87020-020 Maring\'a, PR - Brazil \\
{\small {\em E-mail}:\quad mehernandes@uem.br}
\end{document}